\documentclass[11pt]{article}
\usepackage[russian]{babel}
\usepackage[cp1251]{inputenc}
\usepackage{amsmath,latexsym}
\usepackage[psamsfonts]{amssymb}
\usepackage[dvips]{graphicx}
\usepackage[mathcal]{euscript}
\begin{document}

MSC2020.35G31

\begin{center}
\textbf{ INITIAL-BOUNDARY PROBLEM FOR DEGENERATE
HIGH ORDER EQUATION WITH FRACTIONAL DERIVATIVE
}\\
\textbf{B.Yu.Irgashev}\\
\emph{Namangan Engineering Construction Institute,}\\
\emph{Institute of mathematics Namangan refional department Uzbekistan Academy of Sciences}\\
\emph{Namangan city, Uzbekistan}\\
  \emph{E-mail: bahromirgasev@gmail.com}
\end{center}

\textbf{Annotation..} \emph{The mixed problem for a degenerate high order equation with a fractional derivative in a rectangular domain is considered in the article. The existence of a solution and its uniqueness are shown by the spectral method.}


\begin{center}
\textbf{1.Introduction
}\end{center}

In the region   $D = {D_x} \times {D_y},\,\,{D_x} = \left\{ {x:0\, < x < 1} \right\},\,\,{D_y} = \left\{ {y:\,0 < y < 1} \right\},$ we consider the equation
\[{\left( { - 1} \right)^{k + 1}}D_{0x}^\alpha u\left( {x,y} \right) - {y^m}\frac{{{\partial ^{2k}}u\left( {x,y} \right)}}{{\partial {y^{2k}}}} = 0,\,\,1 < \alpha  < 2,\,\,\,0 \le m < k,\,m \notin N.\eqno(1)\]
where  $ k \in N, D_{0x}^\alpha  - $  is the Riemann-Liouville fractional differentiation operator of order $\alpha$
\[D_{0x}^\alpha u\left( {x,y} \right) = \frac{1}{{\Gamma \left( {2 - \alpha } \right)}}\frac{{{d^2}}}{{d{x^2}}}\int\limits_0^x {\frac{{u\left( {\tau ,y} \right)d\tau }}{{{{\left( {x - \tau } \right)}^{\alpha  - 1}}}}\,} .\]
For equation (1), we consider the problem.

	\textbf{Problem A.}  Find a solution to equation (1) with the conditions:
\[D_{0x}^\alpha u\left( {x,y} \right) \in C\left( {{D_x} \times {{\overline D }_y}} \right),{x^{2 - \alpha }}u\left( {x,y} \right) \in {C^1}\left( {\overline D } \right),\,\frac{{{\partial ^{2k}}u}}{{\partial {y^{2k}}}} \in C\left( {{D_x} \times {{\overline D }_y}} \right),\eqno(2)\]
\[\frac{{{\partial ^{2s}}u\left( {x,0} \right)}}{{\partial {y^{2s}}}} = \frac{{{\partial ^{2s}}u\left( {x,1} \right)}}{{\partial {y^{2s}}}} = 0,\,\,\,0 < x \le 1,\,s = 0,1,...,k - 1,\eqno(3)\]
\[\mathop {\lim }\limits_{x \to 0} {x^{2 - \alpha }}u\left( {x,y} \right) = \psi \left( y \right),\eqno(4)\]
\[\mathop {\lim }\limits_{x \to 0} \frac{\partial }{{\partial x}}\left( {{x^{2 - \alpha }}u\left( {x,y} \right)} \right) = \varphi \left( y \right),\eqno(5)\]
here the functions  $\varphi \left( y \right),\psi \left( y \right) - $  are quite smooth and the natural matching conditions are satisfied for them.\\
Fractional differential equations arise in mathematical modeling of various physical processes and phenomena [1]. Second-order equations, without degeneracy, with partial derivatives of fractional order, were studied in [1] - [5] and others. In these papers, the Cauchy problem, first, second and mixed boundary value problems were considered, a fundamental solution was found, a general representation of the solutions was constructed. Mixed equations and high-order equations with a fractional derivative were studied in [6] - [9]. Fractional order equations with degeneracy were studied in [1, 10].
The study will be carried out by the Fourier method. Earlier, by the Fourier method, boundary value problems for equations with a fractional derivative were studied in [5] - [6], [11].

\begin{center}
\textbf{2. Solution existence
}\end{center}

We are looking for a solution in the form
\[u\left( {x,y} \right) = X\left( x \right)Y\left( y \right).\]
Then with respect to the variable  $y$, taking into account condition (3), we obtain the following spectral problem:
\[\left\{ \begin{array}{l}
{Y^{\left( {2k} \right)}}(y) = {\left( { - 1} \right)^k}\lambda {y^{ - m}}Y(y),\\
{Y^{\left( s \right)}}(0) = {Y^{\left( s \right)}}\left( 1 \right) = 0,\,s = 0,1,...,k - 1.
\end{array} \right.\eqno(6)\]
First, we show that $\lambda  = 0$  is not an eigenvalue. Indeed, consider problem (6) with $\lambda  = 0$
\[\left\{ \begin{array}{l}
{Y^{\left( {2k} \right)}}(y) = 0,\\
{Y^{\left( s \right)}}(0) = {Y^{\left( s \right)}}\left( 1 \right) = 0,\,s = 0,1,...,k - 1.
\end{array} \right.\]
A solution to this problem satisfying the condition ${Y^{\left( s \right)}}(0) = 0,\,s = 0,1,...,k - 1,$  has the form
\[Y\left( y \right) = {c_k}{x^k} + {c_{k + 1}}{x^{k + 1}} + ... + {c_{2k - 1}}{x^{2k - 1}},\]
to determine the unknowns ${c_j},\,\,j = k,k + 1,...,2k - \,1,$  we obtain the system of equations
\[\left\{ \begin{array}{l}
{c_k} + {c_{k + 1}} + ... + {c_{2k - 1}} = 0,\\
k{c_k} + \left( {k + 1} \right){c_{k + 1}} + ... + \left( {2k - 1} \right){c_{2k - 1}} = 0,\\
...\\
k\left( {k - 1} \right)...2 \cdot {c_k} + \left( {k + 1} \right)k...3 \cdot {c_{k + 1}} + ... + \left( {2k - 1} \right)\left( {2k - 2} \right)...\left( {k + 1} \right){c_{2k - 1}} = 0,
\end{array} \right.\]
the main determinant of this system $\Delta $ , has the form
\[\Delta  = \left| {\begin{array}{*{20}{c}}
1&1&{...}&1\\
k&{k + 1}&{...}&{2k - 1}\\
.&.&.&.\\
{k\left( {k - 1} \right)...2}&{\left( {k + 1} \right)k...3}&{...}&{\left( {2k - 1} \right)\left( {2k - 2} \right)...\left( {k + 1} \right)}
\end{array}} \right| = \]
\[ = \prod\limits_{j = 1,j > i}^{k - 1} {\left( {j - i} \right)}  \ne 0.\]
Hence ${c_j} = 0,\,\,j = k,k + 1,...,2k - \,1,$ from here
\[Y\left( y \right) = {c_k}{x^k} + {c_{k + 1}}{x^{k + 1}} + ... + {c_{2k - 1}}{x^{2k - 1}} \equiv 0.\]
We proceed to find the general solution of equation (6) for $\lambda  \ne 0.$ We make a change of variables
\[t = {y^a},\]
\[D_y^{2k}Y\left( y \right) = D_y^{2k - 1}\left( {D_t^1Y\left( y \right)D_y^1t} \right) = \sum\limits_{{k_1} = 0}^{2k - 1} {\left( {C_{2k - 1}^{{k_1}}D_y^{2k - 1 - {k_1}}\left( {D_y^1t} \right)D_y^{{k_1}}\left( {D_t^1Y\left( y \right)} \right)} \right)}  = \]
\[ = D_y^{2k - 1}\left( {D_y^1t} \right)D_t^1Y + \sum\limits_{{k_1} = 1}^{2k - 1} {C_{2k - 1}^{{k_1}}D_y^{2k - 1 - {k_1}}\left( {D_y^1t} \right)D_y^{{k_1} - 1}\left( {D_t^2YD_y^1t} \right)} .\]
 We introduce the following notation: $\overline {{{\left( a \right)}_{j + 1}}}  = a\left( {a - 1} \right)\left( {a - 2} \right)...\left( {a - j} \right)$ , and we assume that $\overline {{{\left( a \right)}_0}}  = 1,$   then $D_y^j\left( {D_y^1t} \right) = \frac{{{{\left( a \right)}_{j + 1}}t}}{{{y^{j + 1}}}}$ . From here we have
\[D_y^{2k}Y\left( y \right) = \frac{{{{\left( a \right)}_{2k}}t}}{{{y^{2k}}}}D_t^1Y\left( y \right) + \sum\limits_{{k_1} = 1}^{2k - 1} {C_{2k - 1}^{{k_1}}\frac{{{{\left( a \right)}_{2k - {k_1}}}}}{{{y^{2k - {k_1}}}}}t} \sum\limits_{{k_2} = 0}^{{k_1} - 1} {C_{{k_1} - 1}^{{k_2}}} D_y^{{k_1} - {k_2} - 1}\left( {D_y^1t} \right)D_y^{{k_2}}\left( {D_t^2Y\left( y \right)} \right) = \]
\[ = \frac{{{{\left( a \right)}_{2k}}}}{{{y^{2k}}}}tD_t^1Y\left( y \right) + \sum\limits_{{k_1} = 1}^{2k - 1} {C_{2k - 1}^{{k_1}}\frac{{{{\left( a \right)}_{2k - {k_1}}}}}{{{y^{2k - {k_1}}}}}tD_y^{{k_1} - 1}\left( {D_y^1t} \right)D_t^2Y\left( y \right)}  + \]
\[ + \sum\limits_{{k_1} = 1}^{2k - 1} {\sum\limits_{{k_2} = 1}^{{k_1} - 1} {C_{2k - 1}^{{k_1}}C_{{k_1} - 1}^{{k_2}}\frac{{{{\left( a \right)}_{2k - {k_1}}}}}{{{y^{2k - {k_1}}}}}\frac{{{{\left( a \right)}_{{k_1} - {k_2}}}}}{{{y^{{k_1} - {k_2}}}}}{t^2}} } D_y^{{k_2} - 1}\left( {D_t^3Y(y)D_y^1t} \right).\]
Continuing this process, we obtain the following formula:
\[D_y^{2k}Y\left( y \right) = {y^{ - 2k}}\sum\limits_{j = 1}^{2k} {\left( {A_{j - 1}^{2k}\left( a \right){t^j}D_t^jY\left( y \right)} \right)},\eqno(7)\]
where
\[A_j^{2k}\left( a \right) = \sum\limits_{{k_1} = j}^{2k - 1} {\sum\limits_{{k_2} = j - 1}^{{k_1} - 1} {...} \sum\limits_{{k_j} = 1}^{{k_{j - 1}} - 1} {\left( {{{\left( a \right)}_{{k_j}}}\prod\limits_{s = 1}^j {C_{{k_{s - 1}} - 1}^{{k_s}}{{\left( a \right)}_{{k_{s - 1}} - {k_s}}}} } \right)} } ,\]
moreover ${k_0} = 2k,\,\,j = \overline {1,2k - 1} ,\,\,{k_1} > {k_2} > ... > {k_j} \ge 1.$
Further, we assume that $A_i^j\left( a \right) = 0$ , for $i \ge j$. We note some properties of the coefficients $A_i^j\left( a \right)$ , that were established in [11].

\textbf{Lemma.}

1. $A_i^{i + 1} = {a^{i + 1}};$

2. $A_i^j = \sum\limits_{k = i}^{j - 1} {C_{j - 1}^k{{\left( a \right)}_{j - k}}A_{i - 1}^k} ;$

3. $A_i^j = a\left( {\left( {i + 1} \right)A_i^{j - 1} + A_{i - 1}^{j - 1}} \right) - \left( {j - 1} \right)A_i^{j - 1};$

4. $A_0^j = {\left( a \right)_j};$

5. $\sum\limits_{j = 1}^s {{{\left( x \right)}_j}A_{j - 1}^s\left( a \right)}  = {\left( {ax} \right)_s}.$\\
We substitute (7) into equation (6)
\[{y^{ - 2k}}\sum\limits_{s = 1}^{2k} {A_{s - 1}^{2k}\left( a \right)} \,{t^s}{Y^{\left( s \right)}}\left( t \right) = {\left( { - 1} \right)^k}\lambda {y^{ - m}}Y\left( t \right),\]
or
\[\sum\limits_{s = 1}^{2k} {A_{s - 1}^{2k}\left( a \right)} \,{t^s}{Y^{\left( s \right)}}\left( t \right) = {\left( { - 1} \right)^k}\lambda {y^{2k - m}}Y\left( t \right).\]
Let be
\[a = 2k - m,\]
then
\[\sum\limits_{s = 1}^{2k} {A_{s - 1}^{2k}\left( a \right)\,} {t^s}{Y^{\left( s \right)}}\left( t \right) = {\left( { - 1} \right)^k}\lambda \,tY\left( t \right).\]
The solution to this equation will be sought in the form of the following series
\[Y\left( t \right) = \sum\limits_{j = 0}^\infty  {{c_j}{t^{\alpha  + j}},} \]
here  $\alpha  - $ unknown parameter. \\
Substituting into the equation, we have
\[\sum\limits_{s = 1}^{2k} {A_{s - 1}^{2k}\left( a \right){t^s}\sum\limits_{j = 0}^\infty  {{c_j}{{\left( {\overline {\alpha  + j} } \right)}_s}{t^{\alpha  + j - s}}} }  = {\left( { - 1} \right)^k}\lambda t\sum\limits_{j = 0}^\infty  {{c_j}{t^{\alpha  + j}},} \]
\[\sum\limits_{j = 0}^\infty  {{c_j}{t^j}\sum\limits_{s = 1}^{2k} {A_{s - 1}^{2k}\left( {2k - m} \right)} {{\left( {\overline {\alpha  + j} } \right)}_s}}  = {\left( { - 1} \right)^k}\lambda \sum\limits_{j = 0}^\infty  {{c_j}{t^{j + 1}},} \]
property 5 of the lemma implies
\[\sum\limits_{s = 1}^{2k} {A_{s - 1}^{2k}\left( {2k - m} \right)} {\overline {\left( {\alpha  + j} \right)} _s} = {\overline {\left( {\left( {2k - m} \right)\left( {\alpha  + j} \right)} \right)} _{2k}},\]
from here
\[\sum\limits_{j = 0}^\infty  {{c_j}{t^j}{{\overline {\left( {\left( {2k - m} \right)\left( {\alpha  + j} \right)} \right)} }_{2k}}}  = {\left( { - 1} \right)^k}\lambda \sum\limits_{j = 0}^\infty  {{c_j}{t^{j + 1}},} \]
now if
\[\alpha  = \frac{s}{{2k - m}},\,\,\,s = 0,1,...,2k - 1,\]
then
\[{c_0} \ne 0,\,\,{c_1} = \frac{{{{\left( { - 1} \right)}^k}\lambda }}{{{{\left( {2k - m} \right)}^{2k}}\prod\limits_{s = 0}^{2k - 1} {\left( {\alpha  - \frac{s}{{2k - m}} + 1} \right)} }}{c_0},\]
because
\[2k - m \notin N,\]
then
\[\prod\limits_{s = 0}^{2k - 1} {\left( {\alpha  - \frac{s}{{2k - m}} + 1} \right)}  \ne 0,\]
further
\[{c_2} = \frac{{{{\left( { - 1} \right)}^k}\lambda {c_1}}}{{{{\left( {2k - m} \right)}^{2k}}\prod\limits_{s = 0}^{2k - 1} {\left( {\alpha  - \frac{s}{{2k - m}} + 2} \right)} }} = \frac{{{{\left( { - 1} \right)}^{2k}}{\lambda ^2}{c_0}}}{{{{\left( {{{\left( {2k - m} \right)}^{2k}}} \right)}^2}\prod\limits_{s = 0}^{2k - 1} {{{\left( {\alpha  - \frac{s}{{2k - m}} + 1} \right)}_2}} }},\]
\[{c_j} = \frac{{{c_0}}}{{\prod\limits_{s = 0}^{2k - 1} {{{\left( {\alpha  - \frac{s}{{2k - m}} + 1} \right)}_j}} }}{\left( {\frac{{{{\left( { - 1} \right)}^k}\lambda }}{{{{\left( {2k - m} \right)}^{2k}}}}} \right)^j},\]
\[Y\left( t \right) = {t^\alpha }\sum\limits_{j = 0}^\infty  {\frac{{{c_0}}}{{\prod\limits_{s = 0}^{2k - 1} {{{\left( {\alpha  - \frac{s}{{2k - m}} + 1} \right)}_j}} }}{{\left( {\frac{{{{\left( { - 1} \right)}^k}\lambda }}{{{{\left( {2k - m} \right)}^{2k}}}}t} \right)}^j} \Rightarrow } \]

\[{Y_i}\left( y \right) = {y^i}\sum\limits_{j = 0}^\infty  {\frac{{{{\left( {\frac{{{{\left( { - 1} \right)}^k}\lambda {y^{2k - m}}}}{{{{\left( {2k - m} \right)}^{2k}}}}} \right)}^j}}}{{j!\prod\limits_{s = 0,s \ne i}^{2k - 1} {{{\left( {\frac{{i - s}}{{2k - m}} + 1} \right)}_j}} }}}, i=0,1,...,2k-1. \]

So in terms of special functions, we got $ 2k $ pieces of linearly independent solutions
\[{Y_i}\left( y \right) = {y^i} \cdot {}_0{F_{2k - 1}}\left[ {\frac{i}{{2k - m}} + 1,...,\frac{{i - \left( {i - 1} \right)}}{{2k - m}} + 1,\frac{{i - \left( {i + 1} \right)}}{{2k - m}} + 1,...,\frac{{i - \left( {2k - 1} \right)}}{{2k - m}} + 1,\frac{{{{\left( { - 1} \right)}^k}\lambda {y^{2k - m}}}}{{{{\left( {2k - m} \right)}^{2k}}}}} \right],\]
\[i = 0,1,...,2k - 1,\]
where
\[{}_p{F_q}\left[ {\begin{array}{*{20}{c}}
{{a_1},...,{a_p},x}\\
{{b_1},...,{b_q}}
\end{array}} \right] = \sum\limits_{k = 0}^\infty  {\frac{{{{\left( {{a_1}} \right)}_k}...{{\left( {{a_p}} \right)}_k}}}{{{{\left( {{b_1}} \right)}_k}...{{\left( {{b_q}} \right)}_k}}}} \frac{{{x^k}}}{{k!}}\]
- is the generalized hypergeometric function, here
 $${\left( a \right)_k} = a\left( {a + 1} \right)...\left( {a + k - 1} \right)$$ - Pohhammer symbol.

	In particular, for $ k = 1 $ we have (${c_0},...,{c_3} - const$)
\[{Y_0}\left( t \right) = {c_0}{\left( {\frac{{\sqrt \lambda  {y^{\frac{{2 - m}}{2}}}}}{{2 - m}}} \right)^{\frac{1}{{2 - m}}}}\sum\limits_{j = 0}^\infty  {\frac{{{{\left( { - 1} \right)}^j}{{\left( {\frac{{2\sqrt \lambda  {y^{\frac{{2 - m}}{2}}}}}{{2\left( {2 - m} \right)}}} \right)}^{2j - \frac{1}{{2 - m}}}}}}{{j!\Gamma \left( {j - \frac{1}{{2 - m}} + 1} \right)}}}  = {c_1}\sqrt y {J_{ - \frac{1}{{2 - m}}}}\left( {\frac{{2\sqrt \lambda  {y^{\frac{{2 - m}}{2}}}}}{{2 - m}}} \right),\]

\[{Y_1}\left( y \right) = {c_2}y\sum\limits_{j = 0}^\infty  {\frac{{{{\left( { - 1} \right)}^j}{{\left( {\frac{{\lambda {y^{2 - m}}}}{{{{\left( {2 - m} \right)}^2}}}} \right)}^j}}}{{j!\Gamma \left( {j + \frac{1}{{2 - m}} + 1} \right)}}}  = {c_3}\sqrt y \sum\limits_{j = 0}^\infty  {\frac{{{{\left( { - 1} \right)}^j}{{\left( {\frac{{2\sqrt \lambda  {y^{\frac{{2 - m}}{2}}}}}{{2\left( {2 - m} \right)}}} \right)}^{2j + \frac{1}{{2 - m}}}}}}{{j!\Gamma \left( {j + \frac{1}{{2 - m}} + 1} \right)}}}  = \]
\[ = {c_3}\sqrt y {J_{\frac{1}{{2 - m}}}}\left( {\frac{{2\sqrt \lambda  {y^{\frac{{2 - m}}{2}}}}}{{2 - m}}} \right),\]
where
\[{J_\nu }\left( z \right) = \sum\limits_{j = 0}^\infty  {\frac{{{{\left( { - 1} \right)}^j}{{\left( {\frac{z}{2}} \right)}^{2j + \nu }}}}{{j!\Gamma \left( {j + \nu  + 1} \right)}}}\]
- Bessel functions [12].\\
Satisfying the boundary conditions, we obtain the condition for the existence of eigenvalues
\[{J_{\frac{1}{{2 - m}}}}\left( {\frac{{2\sqrt \lambda  }}{{2 - m}}} \right) = 0.\]

    Let's get back to the general case. Because  $\left( {2k - m} \right) \notin N$ , then the system of functions $\left\{ {{Y_i}\left( y \right)} \right\}_{i = 0}^{i = 2k - 1}$  - forms a fundamental system of solutions. Hence the general solution of equation (6) has the form
\[Y\left( y \right) = {c_0}{Y_0}\left( y \right) + {c_1}{Y_1}\left( y \right) + ... + {c_{2k - 1}}{Y_{2k - 1}}\left( y \right),\]
from the first boundary condition we obtain
\[Y\left( y \right) = {c_k}{Y_k}\left( y \right) + {c_{k + 1}}{Y_{k + 1}}\left( y \right) + ... + {c_{2k - 1}}{Y_{2k - 1}}\left( y \right),\]
mean
\[Y\left( y \right) = O\left( {{y^k}} \right),\,y \to  + 0,\]
from the second boundary condition we have the system
\[\left\{ \begin{array}{l}
{c_k}{Y_k}\left( 1 \right) + {c_{k + 1}}{Y_{k + 1}}\left( 1 \right) + ... + {c_{2k - 2}}{Y_{2k - 2}}\left( 1 \right) + {c_{2k - 1}}{Y_{2k - 1}}\left( 1 \right) = 0,\\
...\\
\left( {{c_k}{Y_k}\left( y \right) + {c_{k + 1}}{Y_{k + 1}}\left( y \right) + ... + {c_{2k - 2}}{Y_{2k - 2}}\left( y \right) + {c_{2k - 1}}{Y_{2k - 1}}\left( y \right)} \right)_{y = 1}^{\left( {k - 1} \right)} = 0,
\end{array} \right.\]
          equating to zero the main determinant of the system, one can find the eigenvalues of problem (6).  But in view of the complexity of this process, we will proceed in a different way, namely: we reduce problem (6) to the integral equation using the Green function and obtain the necessary estimates for the eigenfunctions. But first, we show that $\lambda  > 0.$ Indeed, we have
\[\int\limits_0^1 {Y\left( y \right){Y^{\left( {2k} \right)}}(y)dy = {{\left( { - 1} \right)}^k}\lambda \int\limits_0^1 {{y^{ - m}}{Y^2}(y)dy} ,} \]
\[\int\limits_0^1 {{{\left( {Y_{}^{\left( k \right)}} \right)}^2}dy = \lambda \int\limits_0^1 {{y^{ - m}}Y_{}^2(y)dy} ,} \]
 because $\lambda  = 0$ is not an eigenvalue, it follows that $\lambda  > 0$. It remains to show the existence of eigenvalues and eigenfunctions of problem (6). The integral equation equivalent to problem (6) has the form
$$Y\left( y \right) = {\left( { - 1} \right)^k}\lambda \int\limits_0^1 {{\xi ^{ - m}}G\left( {y,\xi } \right)Y\left( \xi  \right)d\xi },\eqno(8)$$
where
$$G\left( {y,\xi } \right) =  - \frac{1}{{\left( {2k - 1} \right)!}}\left\{ \begin{array}{l}
{G_1}\left( {y,\xi } \right),\,\,\,\,0 \le y \le \xi ,\\
{G_2}\left( {y,\xi } \right),\,\,\,\,\xi  \le y \le 1,
\end{array} \right.$$
- the Green function of problem (6) (see [13]), here
$${G_1}\left( {y,\xi } \right) = {\left( {1 - \xi } \right)^k}{y^k}\sum\limits_{i = 0}^{k - 1} {\sum\limits_{j = 0}^{k - i - 1} {{{\left( { - 1} \right)}^i}C_{2k - 1}^i} } C_{k - 1 + j}^j{y^{k - i - 1}}{\xi ^{j + i}},$$
$${G_2}\left( {y,\xi } \right) = {\left( {1 - y} \right)^k}{\xi ^k}\sum\limits_{i = 0}^{k - 1} {\sum\limits_{j = 0}^{k - i - 1} {{{\left( { - 1} \right)}^i}C_{2k - 1}^i} } C_{k - 1 + j}^j{\xi ^{k - i - 1}}{y^{j + i}},$$
$$C_n^k = \frac{{n!}}{{k!\left( {n - k} \right)!}}.$$
Rewrite (8) as
$${y^{ - \frac{m}{2}}}Y\left( y \right) = \lambda \int\limits_0^1 {{\xi ^{ - \frac{m}{2}}}\left[ {{{\left( { - 1} \right)}^k}G\left( {y,\xi } \right)} \right]{y^{ - \frac{m}{2}}}\left( {{\xi ^{ - \frac{m}{2}}}Y\left( \xi  \right)} \right)d\xi } ,$$
we introduce the notation
\[\overline Y \left( y \right) = {y^{ - \frac{m}{2}}}Y\left( y \right),\]
$$\overline G \left( {y,\xi } \right) = {\xi ^{ - \frac{m}{2}}}\left[ {{{\left( { - 1} \right)}^k}G\left( {y,\xi } \right)} \right]{y^{ - \frac{m}{2}}},$$
then we have
$$\overline Y \left( y \right) = \lambda \int\limits_0^1 {\overline G \left( {y,\xi } \right)\overline Y \left( \xi  \right)d\xi },\eqno(9)$$
(9) - there is an integral equation with a continuous, in both variables, and a symmetric kernel. According to the theory of equations with symmetric kernels, equation (9) has no more than a countable number of eigenvalues and eigenfunctions. So, problem (6) has eigenvalues ${\lambda _n} > 0,\,\,n = 1,2,...$, and the corresponding eigenfunctions are ${Y_n}(y)$ . Further, we assume that
\[{\left\| {{Y_n}(y)} \right\|^2} = \int\limits_0^1 {{y^{ - m}}Y_n^2(y)dy}  = 1,\]
then, taking into account (9), we have the Bessel inequality
$$\sum\limits_{n = 0}^\infty  {{{\left( {\frac{{{Y_n}\left( y \right)}}{{{\lambda _n}}}} \right)}^2}}  \le \int\limits_0^1 {{y^{ - m}}{G^2}\left( {y,\xi } \right)dy}  < \infty.\eqno(10)$$

Now we find the conditions under which the given function  $\varphi \left( y \right)$  is expanded in a series according to the eigenfunctions ${Y_n}(y)$ . For this we use the Hilbert-Schmidt theorem.

\textbf{Theorem 1.} Let the function $\varphi \left( y \right)$ satisfy the conditions

1. $\varphi \left( y \right) \in {C^{2k}}\left[ {0,1} \right];$

2. ${\varphi ^{\left( i \right)}}\left( 0 \right) = {\varphi ^{\left( i \right)}}\left( 1 \right) = 0,\,i = 0,1,...,k - 1.$\\
Then it can be expanded in a uniformly and absolutely converging series of the form
\[\varphi \left( y \right) = \sum\limits_{n = 1}^\infty  {{\varphi _n}{Y_n}\left( y \right)} ,\]
where
\[{\varphi _n} = \int\limits_0^1 {{y^{ - m}}\varphi \left( y \right){Y_n}(y)dy} .\]
\textbf{Proof.} We show the equality
\[{y^{ - \frac{m}{2}}}\varphi \left( y \right) = \int\limits_0^1 {\overline G \left( {y,\xi } \right)\left( {{{\left( { - 1} \right)}^k}{\xi ^{\frac{m}{2}}}\frac{{{d^{2k}}\varphi \left( \xi  \right)}}{{d{\xi ^{2k}}}}} \right)d\xi } ,\]
really
\[\begin{array}{l}
\int\limits_0^1 {{\xi ^{ - \frac{m}{2}}}\left[ {{{\left( { - 1} \right)}^k}G\left( {y,\xi } \right)} \right]{y^{ - \frac{m}{2}}}\left( {{{\left( { - 1} \right)}^k}{\xi ^{\frac{m}{2}}}\frac{{{d^{2k}}\varphi \left( \xi  \right)}}{{d{\xi ^{2k}}}}} \right)d\xi }  = \\
 = {y^{ - \frac{m}{2}}}\int\limits_0^1 {G\left( {y,\xi } \right)} \frac{{{d^{2k}}\varphi \left( \xi  \right)}}{{d{\xi ^{2k}}}}d\xi  = {y^{ - \frac{m}{2}}}\varphi \left( y \right).
\end{array}\]
Those for the function ${y^{ - \frac{m}{2}}}\varphi \left( y \right)$ the conditions of the Hilbert-Schmidt theorem are satisfied and therefore
\[{y^{ - \frac{m}{2}}}\varphi \left( y \right) = \sum\limits_{n = 1}^\infty  {{y^{ - \frac{m}{2}}}{\varphi _n}{Y_n}\left( y \right)} ,\]
dividing by  ${y^{ - \frac{m}{2}}}$, we have
\[\varphi \left( y \right) = \sum\limits_{n = 1}^\infty  {{\varphi _n}{Y_n}\left( y \right)}.\]
\textbf{Theorem 1 is proved.}

We proceed to solve the equation in the variable $ x $. Taking into account conditions (4), (5), we obtain the following initial problem:
\[\left\{ \begin{array}{l}
D_{0x}^\alpha {X_n}\left( x \right) =  - {\lambda _n}{X_n}\left( x \right),\\
\mathop {\lim }\limits_{x \to 0} \left( {{x^{2 - \alpha }}{X_n}\left( x \right)} \right) = {\psi _n},\\
\mathop {\lim }\limits_{x \to 0} \frac{d}{{dx}}\left( {{x^{2 - \alpha }}{X_n}\left( x \right)} \right) = {\varphi _n}.
\end{array} \right.\eqno(11)\]
where\\
$${\psi _n} = \int\limits_0^1 {\psi \left( y \right){y^{ - m}}{Y_n}\left( y \right)dy\,\,} ,$$
$${\varphi _n} = \int\limits_0^1 {\varphi \left( y \right){y^{ - m}}{Y_n}\left( y \right)dy\,\,} .$$
We will search for solution (11) as a series
\[{X_n}\left( x \right) = \sum\limits_{j = 0}^\infty  {{c_j}{x^{\gamma j + \beta }}} ,\]
where ${c_j},\beta ,\gamma- $ are still unknown real numbers. \\
Formally, we have (the legality of rearranging the series and the integral will follow below)
\[D_{0x}^\alpha {X_n}\left( {x} \right) = \frac{1}{{\Gamma \left( {2 - \alpha } \right)}}\frac{{{d^2}}}{{d{x^2}}}\int\limits_0^x {\frac{{\sum\limits_{j = 0}^\infty  {{c_j}{\tau ^{\gamma j + \beta }}} d\tau }}{{{{\left( {x - \tau } \right)}^{\alpha  - 1}}}}\, = } \frac{1}{{\Gamma \left( {2 - \alpha } \right)}}\sum\limits_{j = 0}^\infty  {{c_j}} \frac{{{d^2}}}{{d{x^2}}}\int\limits_0^x {\frac{{{\tau ^{\gamma j + \beta }}d\tau }}{{{{\left( {x - \tau } \right)}^{\alpha  - 1}}}}\, = } \]\[ = \frac{1}{{\Gamma \left( {2 - \alpha } \right)}}\sum\limits_{j = 0}^\infty  {{c_j}} \frac{{{d^2}}}{{d{x^2}}}\int\limits_0^1 {\frac{{{x^{\gamma j + \beta  + 1}}{z^{\gamma j + \beta }}dz}}{{{x^{\alpha  - 1}}{{\left( {1 - z} \right)}^{\alpha  - 1}}}}\, = } \frac{1}{{\Gamma \left( {2 - \alpha } \right)}}\sum\limits_{j = 0}^\infty  {{c_j}} \frac{{{d^2}}}{{d{x^2}}}{x^{\gamma j + \beta  + 2 - \alpha }}\int\limits_0^1 {\frac{{{z^{\gamma j + \beta }}dz}}{{{{\left( {1 - z} \right)}^{\alpha  - 1}}}}\, = } \]
\[ = \sum\limits_{j = 0}^\infty  {{c_j}{{\overline {\left( {\gamma j + \beta  + 2 - \alpha } \right)} }_2}} {x^{\gamma j + \beta  - \alpha }}\frac{{\Gamma \left( {\gamma j + \beta  + 1} \right)}}{{\Gamma \left( {\gamma j + \beta  + 3 - \alpha } \right)}}.\]
Substituting the last expression in (11), we obtain
\[\begin{array}{l}
\sum\limits_{j = 0}^\infty  {{c_j}{{\overline {\left( {\gamma j + \beta  + 2 - \alpha } \right)} }_2}} {x^{\gamma j + \beta  - \alpha }}\frac{{\Gamma \left( {\gamma j + \beta  + 1} \right)}}{{\Gamma \left( {\gamma j + \beta  + 3 - \alpha } \right)}} = \\
 =  - {\lambda _n}\sum\limits_{j = 0}^\infty  {{c_j}{x^{\gamma j + \beta }}}  \Rightarrow
\end{array}\]
\[\sum\limits_{j = 0}^\infty  {{c_j}{{\overline {\left( {\gamma j + \beta  + 2 - \alpha } \right)} }_2}} {x^{\gamma \left( {j - 1} \right) + \gamma  - \alpha }}\frac{{\Gamma \left( {\gamma j + \beta  + 1} \right)}}{{\Gamma \left( {\gamma j + \beta  + 3 - \alpha } \right)}} = \]
\[ =  - {\lambda _n}\sum\limits_{j = 0}^\infty  {{c_j}{x^{\gamma j}}}  \Rightarrow \gamma  - \alpha  = 0 \Rightarrow \gamma  = \alpha  \Rightarrow \]
\[\sum\limits_{j = 0}^\infty  {{c_j}{{\overline {\left( {\alpha j + \beta  + 2 - \alpha } \right)} }_2}} {x^{\alpha \left( {j - 1} \right)}}\frac{{\Gamma \left( {\gamma j + \beta  + 1} \right)}}{{\Gamma \left( {\gamma j + \beta  + 3 - \alpha } \right)}} = \]
\[ =  - {\lambda _n}\sum\limits_{j = 0}^\infty  {{c_j}{x^{\alpha j}}}  \Rightarrow \,\,{\beta _1} = \alpha  - 1,\,{\beta _2} = \alpha  - 2.\]
Let be ${\beta _1} = \alpha  - 1$
\[\sum\limits_{j = 1}^\infty  {{c_j}{{\overline {\left( {\alpha j + 1} \right)} }_2}} {x^{\alpha \left( {j - 1} \right)}}\frac{{\Gamma \left( {\alpha j + \alpha } \right)}}{{\Gamma \left( {\alpha j + 2} \right)}} =  - {\lambda _n}\sum\limits_{j = 0}^\infty  {{c_j}{x^{\alpha j}}}  \Rightarrow \]
\[{c_j} = \frac{{ - {\lambda _n}\left( {\alpha j} \right)\left( {\alpha j + 1} \right)\Gamma \left( {\alpha j} \right)}}{{\left( {\alpha j} \right)\left( {\alpha j + 1} \right)\Gamma \left( {\alpha \left( {j + 1} \right)} \right)}}{c_{j - 1}} = \frac{{{{\left( { - {\lambda _n}} \right)}^j}\Gamma \left( \alpha  \right)}}{{\Gamma \left( {\alpha \left( {j + 1} \right)} \right)}}{c_0} \Rightarrow \]
The first solution would be to
\[{X_{1n}}\left( x \right) = {x^{\alpha  - 1}}\sum\limits_{j = 0}^\infty  {\frac{{{{\left( { - {\lambda _n}{x^\alpha }} \right)}^j}}}{{\Gamma \left( {\alpha j + \alpha } \right)}}.} \]
This series converges absolutely and uniformly for fixed values of ${\lambda _n}$ and for limited values of $ x $.Indeed, on the basis of d'Alembert, we have
\[\mathop {\lim }\limits_{j \to  + \infty } \left( {\frac{{{{\left( {  {\lambda _n}{x^\alpha }} \right)}^{j + 1}}}}{{\Gamma \left( {\alpha j + 2\alpha } \right)}}:\frac{{{{\left( {  {\lambda _n}{x^\alpha }} \right)}^j}}}{{\Gamma \left( {\alpha j + \alpha } \right)}}} \right) = \left( {  {\lambda _n}{x^\alpha }} \right)\mathop {\lim }\limits_{j \to  + \infty } \frac{{\Gamma \left( {\alpha j + \alpha } \right)}}{{\Gamma \left( {\alpha j + 2\alpha } \right)}} = \]
\[ = \left( { {\lambda _n}{x^\alpha }} \right)\mathop {\lim }\limits_{j \to  + \infty } O{\left( {\alpha j} \right)^{\alpha  - 2\alpha }} = 0,\]
here we used the relation from [14]
\[\frac{{\Gamma \left( {z + \alpha } \right)}}{{\Gamma \left( {z + \beta } \right)}} = {z^{\alpha  - \beta }}\left[ {1 + \frac{{\left( {\alpha  - \beta } \right)\left( {\alpha  - \beta  - 1} \right)}}{{2z}} + O\left( {{z^{ - 2}}} \right)} \right].\]
So the permutation of the series and the integral in the above was legal.\\
Now let ${\beta _2} = \alpha  - 2$, then
\[\begin{array}{l}
\sum\limits_{j = 0}^\infty  {{c_j}{{\overline {\left( {\alpha j} \right)} }_2}} {x^{\alpha \left( {j - 1} \right)}}\frac{{\Gamma \left( {\alpha j + \alpha  - 1} \right)}}{{\Gamma \left( {\alpha j + 1} \right)}} = \\
 =  - {\lambda _n}\sum\limits_{j = 0}^\infty  {{c_j}{x^{\alpha j}}}  \Rightarrow
\end{array}\]
\[\sum\limits_{j = 1}^\infty  {{{\overline {\left( {\alpha j} \right)} }_2}} {x^{\alpha \left( {j - 1} \right)}}\frac{{\Gamma \left( {\alpha j + \alpha  - 1} \right)}}{{\Gamma \left( {\alpha j + 1} \right)}} =  - {\lambda _n}\sum\limits_{j = 0}^\infty  {{c_j}{x^{\alpha j}}}  \Rightarrow \]
\[{c_j} = \frac{{ - {\lambda _n}\Gamma \left( {\alpha j + 1} \right)}}{{\left( {\alpha j} \right)\left( {\alpha j - 1} \right)\Gamma \left( {\alpha \left( {j + 1} \right) - 1} \right)}}{c_{j - 1}} = \frac{{ - {\lambda _n}\Gamma \left( {\alpha j - 1} \right)}}{{\Gamma \left( {\alpha \left( {j + 1} \right) - 1} \right)}}{c_{j - 1}} = \]
\[ = ... = \frac{{{{\left( { - {\lambda _n}} \right)}^j}\Gamma \left( {\alpha  - 1} \right)}}{{\Gamma \left( {\alpha \left( {j + 1} \right) - 1} \right)}}{c_0} \Rightarrow \]
the second linearly independent solution will be
\[{X_{2n}}\left( x \right) = {x^{\alpha  - 2}}\sum\limits_{j = 0}^\infty  {\frac{{{{\left( { - {\lambda _n}{x^\alpha }} \right)}^j}}}{{\Gamma \left( {\alpha j + \alpha  - 1} \right)}}.} \]
So, the general solution of equation (11) has the form
\[{X_n}\left( x \right) = {d_1}{x^{\alpha  - 1}}\sum\limits_{j = 0}^\infty  {\frac{{{{\left( { - {\lambda _n}{x^\alpha }} \right)}^j}}}{{\Gamma \left( {\alpha j + \alpha } \right)}} + {d_2}} {x^{\alpha  - 2}}\sum\limits_{j = 0}^\infty  {\frac{{{{\left( { - {\lambda _n}{x^\alpha }} \right)}^j}}}{{\Gamma \left( {\alpha j + \alpha  - 1} \right)}},} \,\,{d_1},{d_2} - const,\]
or in terms of special functions
\[{X_n}\left( x \right) = {d_1}{x^{\alpha  - 1}}{E_{{1 \mathord{\left/
 {\vphantom {1 \alpha }} \right.
 \kern-\nulldelimiterspace} \alpha }}}\left( { - {\lambda _n}{x^\alpha },\alpha } \right) + {d_2}{x^{\alpha  - 2}}{E_{{1 \mathord{\left/
 {\vphantom {1 \alpha }} \right.
 \kern-\nulldelimiterspace} \alpha }}}\left( { - {\lambda _n}{x^\alpha },\alpha  - 1} \right),\eqno(12)\]
 where
 \[{E_{{1 \mathord{\left/
 {\vphantom {1 \alpha }} \right.
 \kern-\nulldelimiterspace} \alpha }}}\left( {z,\mu } \right) = \sum\limits_{j = 0}^\infty  {\frac{{{z^j}}}{{\Gamma \left( {\alpha j + \mu } \right)}}} \]
 - 	Mittag-Leffler function [1].\\
 Note that the representation in the form (12) coincides with the known results (see [1] and others) obtained using the properties of the Mittag-Leffler function.\\
 Satisfying the initial conditions, we obtain a solution to problem (11) in the form
 \[{X_n}\left( x \right) = \Gamma \left( \alpha  \right){\varphi _n}{x^{\alpha  - 1}}{E_{{1 \mathord{\left/
 {\vphantom {1 \alpha }} \right.
 \kern-\nulldelimiterspace} \alpha }}}\left( { - {\lambda _n}{x^\alpha },\alpha } \right) + \Gamma \left( {\alpha  - 1} \right){\psi _n}{x^{\alpha  - 2}}{E_{{1 \mathord{\left/
 {\vphantom {1 \alpha }} \right.
 \kern-\nulldelimiterspace} \alpha }}}\left( { - {\lambda _n}{x^\alpha },\alpha  - 1} \right),\]
 this representation implies the uniqueness of the solution to problem (11).\\
 Given the estimate (see [15], p.136)
 \[\left| {{E_{{1 \mathord{\left/
 {\vphantom {1 \alpha }} \right.
 \kern-\nulldelimiterspace} \alpha }}}\left( { - {\lambda _n}{x^\alpha },\alpha } \right)} \right| \le \frac{M}{{1 + {\lambda _n}{x^\alpha }}},\,0 < M - {\mathop{\rm const}\nolimits} ,\,\alpha  < 2,\]
 we have
 \[\left| {{X_n}\left( x \right)} \right| \le {M_1}\left( {\left| {{\varphi _n}} \right| + \frac{{\left| {{\psi _n}} \right|}}{x}} \right),\,0 < {M_1} - const.\]
Thus, the formal solution of the problem A has the form
\[u\left( {x,y} \right) = \sum\limits_{n = 0}^\infty  {{X_n}\left( x \right){Y_n}\left( y \right)},\eqno(13)\]
We show that (13) is a classical solution of equation (1), indeed
\[\left| {D_{0x}^\alpha u\left( {x,y} \right)} \right| \le \sum\limits_{n = 0}^\infty  {\left| {D_{0x}^\alpha {X_n}\left( x \right)} \right|\left| {{Y_n}\left( y \right)} \right|}  = \sum\limits_{n = 0}^\infty  {\left| {{\lambda _n}{X_n}\left( x \right)} \right|\left| {{Y_n}\left( y \right)} \right|}  \le \]
\[ \le {M_1}\left( {\sum\limits_{n = 0}^\infty  {\left| {{\lambda _n}{\varphi _n}} \right|\left| {{Y_n}\left( y \right)} \right|}  + \frac{1}{x}\sum\limits_{n = 0}^\infty  {\left| {{\lambda _n}{\psi _n}} \right|\left| {{Y_n}\left( y \right)} \right|} } \right),\,\,\]
We show the convergence of the first term $\sum\limits_{n = 0}^\infty  {\left| {{\lambda _n}{\varphi _n}} \right|\left| {{Y_n}\left( y \right)} \right|}$ , convergence of the second: $\sum\limits_{n = 0}^\infty  {\left| {{\lambda _n}{\psi _n}} \right|\left| {{Y_n}\left( y \right)} \right|}$ , is shown in the same way. So
 \[\sum\limits_{n = 0}^\infty  {\left| {{\lambda _n}{\varphi _n}} \right|\left| {{Y_n}\left( y \right)} \right|}  = \sum\limits_{n = 1}^\infty  {\left| {\lambda _n^2{\varphi _n}} \right|\left| {\frac{{{Y_n}\left( y \right)}}{{{\lambda _n}}}} \right|}  \le \sqrt {\sum\limits_{n = 1}^\infty  {\left| {\lambda _n^4\varphi _n^2} \right|} } \,\,\sqrt {\sum\limits_{n = 1}^\infty  {\left| {\frac{{Y_n^2\left( y \right)}}{{\lambda _n^2}}} \right|} } ,\]
we have
\[\begin{array}{l}
{\varphi _n} = \int\limits_0^1 {{y^{ - m}}\varphi \left( y \right){Y_n}(y)dy}  = \frac{{{{\left( { - 1} \right)}^k}}}{{{\lambda _n}}}\int\limits_0^1 {\varphi \left( y \right)Y_n^{\left( {2k} \right)}(y)dy}  = \\
 = \frac{{{{\left( { - 1} \right)}^k}}}{{{\lambda _n}}}\int\limits_0^1 {{\varphi ^{\left( {2k} \right)}}\left( y \right)Y_n^{}(y)dy}  =
\end{array}\]
\[ = \frac{1}{{\lambda _n^2}}\int\limits_0^1 {{y^m}{\varphi ^{\left( {2k} \right)}}\left( y \right)Y_n^{\left( {2k} \right)}(y)dy}  = \frac{1}{{\lambda _n^2}}\int\limits_0^1 {{y^m}{{\left( {{y^m}{\varphi ^{\left( {2k} \right)}}\left( y \right)} \right)}^{\left( {2k} \right)}}{Y_n}(y){y^{ - m}}dy} ,\]
Now we apply the Bessel inequality
\[\sum\limits_{n = 1}^\infty  {\left| {\lambda _n^4\varphi _n^2} \right|}  \le \int\limits_0^1 {{{\left( {{{\left( {{y^m}{\varphi ^{\left( {2k} \right)}}\left( y \right)} \right)}^{\left( {2k} \right)}}} \right)}^2}{y^{-m}}dy < \infty .} \]
Now, in order for the calculations made above to be legal, we impose the following restrictions on the function $\varphi \left( y \right)$:
\[{\varphi ^{\left( s \right)}}\left( 0 \right) = {\varphi ^{\left( s \right)}}\left( 1 \right) = 0,\,\,\varphi \left( y \right) \in {C^{2k}}\left[ {0,1} \right],\]
\[{\left( {{y^m}{\varphi ^{\left( {2k} \right)}}\left( y \right)} \right)^{\left( s \right)}}\left( 0 \right) = {\left( {{y^m}{\varphi ^{\left( {2k} \right)}}\left( y \right)} \right)^{\left( s \right)}}\left( 1 \right) = 0,\,\,\]
\[{y^m}{\varphi ^{\left( {2k} \right)}}\left( y \right) \in {C^{2k}}\left[ {0,1} \right],s = 0,1,...,k - 1.\]
Taking into account (10) and (14), we obtain that the series
\[D_{0x}^\alpha u\left( {x,y} \right) = \sum\limits_{n = 0}^\infty  {D_{0x}^\alpha {X_n}\left( x \right){Y_n}\left( y \right)}  - \]
converges uniformly. The uniform convergence of the series is proved in a similar way
\[\frac{{{\partial ^{2k}}u\left( {x,y} \right)}}{{\partial {y^{2k}}}} = \sum\limits_{n = 0}^\infty  {{X_n}\left( x \right)\frac{{{\partial ^{2k}}{Y_n}\left( y \right)}}{{\partial {y^{2k}}}}}  = {\left( { - 1} \right)^k}{y^{ - m}}\sum\limits_{n = 0}^\infty  {{\lambda _n}{X_n}\left( x \right){Y_n}\left( y \right)} .\]
So the following theorem holds.
	
\textbf{Theorem 2.} Let the function  $\tau \left( y \right),$  где $\tau \left( y \right) = \varphi \left( y \right)$   or $\tau \left( y \right) = \psi \left( y \right)$ ,  satisfies the following conditions:
\[\tau \left( y \right) \in {C^{2k}}\left[ {0,1} \right],\,{\tau ^{\left( s \right)}}\left( 0 \right) = {\tau ^{\left( s \right)}}\left( 1 \right) = 0,\]
\[{\left( {{y^m}{\tau ^{\left( {2k} \right)}}\left( y \right)} \right)^{\left( s \right)}}\left( 0 \right) = {\left( {{y^m}{\tau ^{\left( {2k} \right)}}\left( y \right)} \right)^{\left( s \right)}}\left( 1 \right) = 0,\,{y^m}{\tau ^{\left( {2k} \right)}}\left( y \right) \in {C^{2k}}\left[ {0,1} \right],s = 0,1,...,k - 1.\]
Then a solution to Problem A exists.
\begin{center}
\textbf{3. Uniqueness of solution
}\end{center}

 Let the function $u\left( {x,y} \right)$ be a solution to Problem A with zero initial and boundary conditions. We consider its Fourier coefficients with respect to the system of eigenfunctions of problem (6)
\[{u_n}\left( x \right) = \int\limits_0^1 {{y^{ - m}}u\left( {x,y} \right){Y_n}(y)dy} ,\]
it is easy to show that ${u_n}\left( x \right)$ is a solution to the problem
\[\left\{ \begin{array}{l}
D_{0x}^\alpha {u_n}\left( x \right) =  - {\lambda _n}{u_n}\left( x \right),\\
\mathop {\lim }\limits_{x \to 0} \left( {{x^{2 - \alpha }}{u_n}\left( x \right)} \right) = 0,\\
\mathop {\lim }\limits_{x \to 0} \frac{d}{{dx}}\left( {{x^{2 - \alpha }}{u_n}\left( x \right)} \right) = 0.
\end{array} \right.\]
This problem has only a zero solution, i.e.
\[\int\limits_0^1 {{y^{ - m}}u\left( {x,y} \right){Y_n}(y)dy}  = 0,\,\forall n.\]
Because \\

$\overline G \left( {y,\xi } \right) - $ symmetrical, continuous\\
\[\int\limits_0^1 {{{\overline G }^2}\left( {y,\xi } \right)d\xi }  < \infty ,\,\int\limits_0^1 {{{\overline G }^2}\left( {y,\xi } \right)dy}  < \infty ,\,\int\limits_0^1 {\int\limits_0^1 {{{\overline G }^2}\left( {y,\xi } \right)dyd\xi } }  < \infty ,\,{\lambda _n} > 0,\,\forall n,\]
then the conditions of the Mercer theorem are satisfied and
\[\overline G \left( {y,\xi } \right) = \sum\limits_{n = 0}^\infty  {\frac{{\overline {{Y_n}} \left( y \right)\overline {{Y_n}} \left( \xi  \right)}}{{{\lambda _n}}}} .\]
From here we have
\[{y^{ - \frac{m}{2}}}u\left( {x,y} \right) = \int\limits_0^1 {\overline G \left( {y,\xi } \right)\left( {{{\left( { - 1} \right)}^k}{\xi ^{\frac{m}{2}}}\frac{{{\partial ^{2k}}u\left( {x,\xi } \right)}}{{\partial {\xi ^{2k}}}}} \right)d\xi }  = \]
\[ = {\left( { - 1} \right)^k}\int\limits_0^1 {\sum\limits_{n = 0}^\infty  {\frac{{\overline {{Y_n}} \left( y \right)\overline {{Y_n}} \left( \xi  \right)}}{{{\lambda _n}}}} \left( {{\xi ^{\frac{m}{2}}}\frac{{{\partial ^{2k}}u\left( {x,\xi } \right)}}{{\partial {\xi ^{2k}}}}} \right)d\xi }  = \]
\[ = {\left( { - 1} \right)^k}\sum\limits_{n = 0}^\infty  {\frac{{{y^{ - \frac{m}{2}}}{Y_n}\left( y \right)}}{{{\lambda _n}}}} \int\limits_0^1 {{\xi ^{ - \frac{m}{2}}}{Y_n}\left( \xi  \right){\xi ^{\frac{m}{2}}}\frac{{{\partial ^{2k}}u\left( {x,\xi } \right)}}{{\partial {\xi ^{2k}}}}d\xi }  = \]
because since the series converges uniformly, then we can interchange the signs of integration and the sum
\[ = {\left( { - 1} \right)^k}\sum\limits_{n = 0}^\infty  {\frac{{{y^{ - \frac{m}{2}}}{Y_n}\left( y \right)}}{{{\lambda _n}}}} \int\limits_0^1 {{Y_n}\left( \xi  \right)\frac{{{\partial ^{2k}}u\left( {x,\xi } \right)}}{{\partial {\xi ^{2k}}}}d\xi }  = \]
\[ = {\left( { - 1} \right)^k}\sum\limits_{n = 0}^\infty  {\frac{{{y^{ - \frac{m}{2}}}{Y_n}\left( y \right)}}{{{\lambda _n}}}} \int\limits_0^1 {Y_n^{\left( {2k} \right)}\left( \xi  \right)u\left( {x,\xi } \right)d\xi }  = \]
\[ = {\left( { - 1} \right)^k}\sum\limits_{n = 0}^\infty  {\frac{{{y^{ - \frac{m}{2}}}{Y_n}\left( y \right)}}{{{\lambda _n}}}} \int\limits_0^1 {{\lambda _n}{{\left( { - 1} \right)}^k}{\xi ^{ - m}}{Y_n}\left( \xi  \right)u\left( {x,\xi } \right)d\xi }  = \]
\[ = {y^{ - \frac{m}{2}}}\sum\limits_{n = 0}^\infty  {{Y_n}\left( y \right)} \int\limits_0^1 {{\xi ^{ - m}}{Y_n}\left( \xi  \right)u\left( {x,\xi } \right)d\xi }  = 0 \Rightarrow \]
\[u\left( {x,y} \right) \equiv 0.\]

\begin{center}
\textbf{References}
\end{center}

1. A. M. Nakhushev; Drobniye ischisleniya i yego primenenie. [Fractional Calculus and its Applications], Moscow: FIZMATLIT, 2003.

2. A. V. Pskhu; Uravneniya v chastnykh proizvodnykh drobnogo poryadka. [Partial differential equations of fractional order], Moscow: NAUKA, 2005.

3. Mainardi F; The fundamental solutions for the fractional diffusion-wave equation. Appl. Math. Lett. 9 (6) (1996), pp. 23–28.

4. A. V. Pskhu; The fundamental solution of a diffusion-wave equation of fractional order. Izv. Math., 73:2 (2009), pp. 351–392

5. O. P. Agrawal;  Solution for a fractional diffusion-wave equation defined in a bounded domain. Nonlinear Dynam. 29 (1) (2002), pp. 145–155.

6. O.Kh. Masaeva; Uniqueness of solutions to Dirichlet problems for generalized Lavrent’ev – Bitsadze equations with a fractional derivative. Electronic Journal of Differential Equations, Vol. 2017 (2017), No. 74, pp. 1–8.

7. B. J. Kadirkulov;  Boundary problems for mixed parabolic-hyperbolic equations with two lines of changing type and fractional derivative.  Electronic Journal of Differential Equations, Vol. 2014 (2014), No. 57, pp. 1–7.

8. A.V. Pskhu; Electronic Journal of Differential Equations, Vol. 2019 (2019), No. 21, pp. 1-13.

9. A.S. Berdyshev , B.E. Eshmatov , B.J. Kadirkulov ; Boundary value problems for fourth-order mixed type equation with fractional derivative.Electronic Journal of Differential Equations, Vol. 2016 (2016), No. 36, pp. 1–11.

10. А. N.Аrtyushin;  Fractional integral inequalities and their applications to degenerate differential
equations with the caputo fractional derivative.  Siberian Mathematical Journal, 61(2) (2020), pp. 208–221.

11. B. Yu.Irgashev ; On partial solutions of one equation with multiple characteristics and some properties of the fundamental solution. Ukrainian Mathematical Journal. 68(6) (2016), pp. 868-893.

12.H.Bateman , A.Erdelyi ; Bateman Manuscript Project. //Higher Transcendental Functions. – McGraw-Hill New York, 1953. – Vol. 2.

13. Yu. P. Apakov, B.Yu.  Irgashev; Boundary-value problem for a degenerate high-odd-order equation. Ukrainian Mathematical Journal.66(10)(2015), pp 1475-1490.

14. H.Bateman , A.Erdelyi  Bateman Manuscript Project. //Higher Transcendental Functions. – McGraw-Hill New York, 1953. – Vol. 1.

15. M. M. Dzhrbashian; Integral transformation and representation of functions in complex domain, Moscow, 1966, 672 p.

\end{document}